# Integrated design optimization of structural bending filter and gain schedules for rocket attitude control system

Sang-Il Lee[a,b], Jaemyung Ahn[a,*], Woong-Rae Roh[b]

[a] *Department of Aerospace Engineering, Korea Advanced Institute of Science and Technology, 291 Daehak-ro, Daejeon 34141, Republic of Korea*
[b] *KSLV-II R&D Head Office, Korea Aerospace Research Institute, 169-84 Gwahak-ro, Daejoen 34133, Republic of Korea*

## Abstract

This paper proposes an integrated design optimization framework for the gain schedules and bending filter for the longitudinal control of a rocket during its ascent flight. Dynamic models representing the pitch/yaw motion of a rocket considering the elements such as the rigid body dynamics, aerodynamics, sloshing, bending, sensor/actuator, and flight computer are introduced. The linear proportional and differential (PD) control law with scheduled (time-varying) gains and bending filter parameters are identified as key decision variables for stabilizing the pitch/yaw motion of the rocket. The integrated optimal design problem that determines the decision variables to minimize the worst-case peak associated with the first bending mode with constraints on the stability margins during the flight of the rocket is mathematically formulated. A case study on design of gain schedules and bending filter for an actual sounding rocket using the proposed framework is conducted to demonstrate its effectiveness.

## 1. Introduction

The challenges in the pitch and yaw control of a rocket during its ascent flight are well known and have been studied by many researchers to date. The subject involves a number of complications primarily coming from the unstable and complex dynamics of the rocket system. Two important sources of the complications for these challenges are the large variations in the system/environmental parameters (e.g., weight, rotational inertia, dynamic pressure) and the structural bending the rocket body associated with its flexibility and slender shape. Traditionally, these two challenges have been addressed by using the gain scheduling technique and the bending filter, which should secure a certain level of stability margin in the attitude control system throughout the whole flight period.

---

* Corresponding author. E-mail: jaemyung.ahn@kaist.ac.kr

A number of studies on these issues can be found in the literature on attitude control of a rocket. Gain scheduling and bending filter design of early launch vehicles such as Saturn, Delta, and Space Shuttle are presented in Haeussermann (1970), Frosch et al. (1968), Simmons et al. (1973), and Schleich (1982). These techniques were applied to H-II (Mori, 1999) and M-V (Morita & Kawaguchi, 2001) launch vehicles. Clement et al. (2001; 2005) presented a gain scheduling methodology based on the Youla parameterization and optimization with linear matrix inequality (LMI) constraints for Ariane 5 launch vehicle. Choi & Bang (2000) proposed an adaptive notch filter for the attitude control of a sounding rocket and Oh et al. (2008) conducted an experiment to demonstrate its effectiveness. As a relatively recent studies, bending filter design and the gain scheduling Ares-I rocket were discussed by Jang et al. (2008) and Jang et al. (2011), respectively.

The aforementioned studies on the gain scheduling and bending filter design problems addressed the two subjects separately and sequentially. That is, the gain profiles for the attitude control system considering the low-frequency dynamics (rigid body and aerodynamics) of the rocket and parameter variations are designed first, and the design of the structural filter that can stabilize the high-frequency bending modes is conducted. This study introduces a new approach to design the gain schedules and bending filter simultaneously using an integrated optimization framework, which provides a systematic procedure that can prevent unnecessary improve the filtering performance (defined as the attenuation of the first bending mode frequency).

The contribution of this is study is twofold. First, an integrated framework to design/optimize simultaneously the bending filter and the gain schedules of the attitude system for the ascent flight of a rocket was proposed. The framework identifies the key parameters of the bending filter and gain profiles, and specify them as the decision variables of the design. Mathematical formulations for the optimal design problem are developed by defining the objective function and constraints associated with the stability of the attitude control loop of the rocket, which can be solved to optimality by using the nonlinear programming (NLP) technique. Second, a case study on design of the bending filter and gain schedules using the proposed framework for an actual sounding rocket system (Korea sounding rocket III, KSR-III) was presented. Models describing an actual rocket dynamics (e.g., rigid body, aerodynamics, sloshing of the liquid propellant, and structural bending), sensing (e.g., sensor and V/F converter), actuator dynamics, and control law are introduced. The optimized filter and gain schedule designs were compared with those used for actual flight of KSR-III, which were obtained by sequential design of gain schedules and bending filter to demonstrate the effectiveness of the proposed framework.

The remainder of this paper is structured as follows. The models representing the elements of the rocket's attitude control system are introduced in Section 2. The feedback control law for the attitude control system and various bending filter structures are presented in Section 3. Section 4 presents the mathematical formulations describing the integrated design problem to optimize the gain schedules and bending filter parameters. An integrated gain schedules and filter design case study for a sounding rocket is introduced in Section 5. Finally, Section 6 provides the conclusions of this study and discusses potential future work opportunities.

## 2. Modeling the attitude control system of a liquid-propellant rocket

This section introduces the equations of motion for describing the model of an attitude control system of a liquid-propellant rocket. The rocket dynamics model (the rigid body mode, the bending mode, the tail-wags-dog (TWD) effects, and the sloshing are introduced), the actuator model for the thrust vector control (TVC), the sensor (dynamically tuned gyro, DTG) model, and the flight computer model (the voltage-to frequency converter (V/F), the digital-to-analog converter (DAC), and the zero-order holder (ZOH)) are presented. The differential equation and/or the transfer functions for the aforementioned elements are presented, which are used for classical linear stability analysis. The overall structure of the control system is shown in Fig. 1.

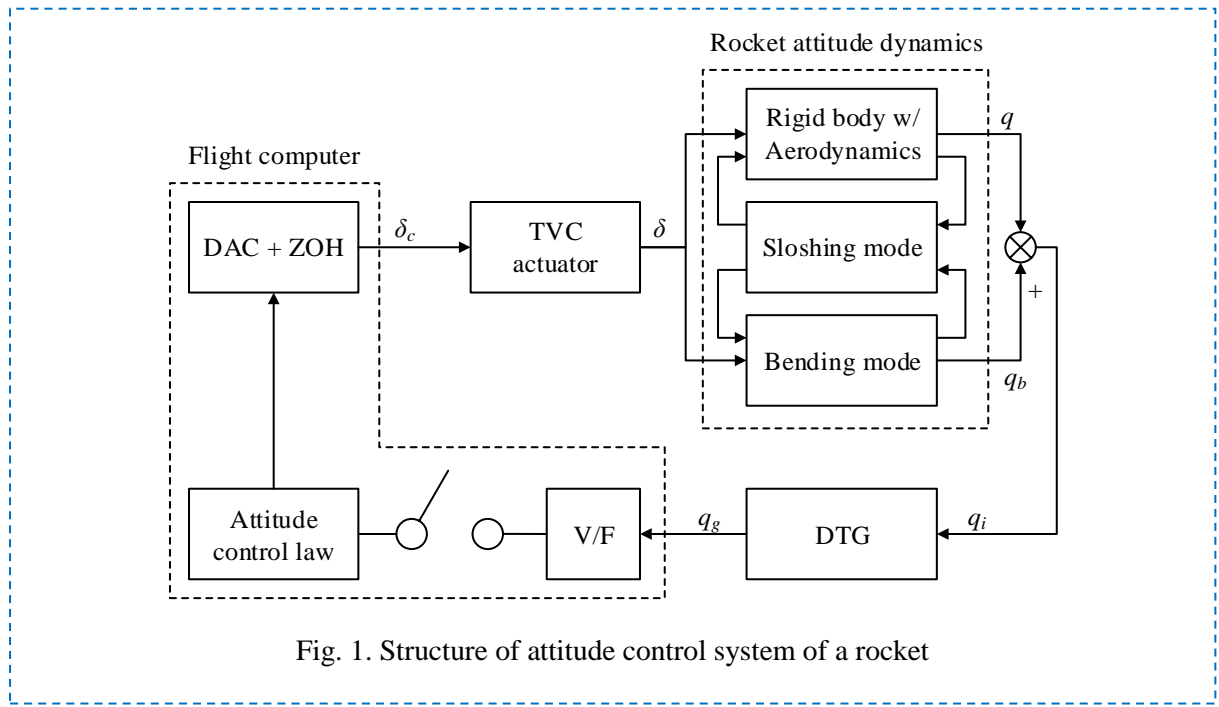

Fig. 1. Structure of attitude control system of a rocket

### 2.1. Attitude dynamics of a rocket

The attitude dynamics of a rocket consists of four parts: the rigid body motion, the bending motion, the TWD mode, and the sloshing. A reference coordinate for the dynamic modeling is shown in Fig. 2.

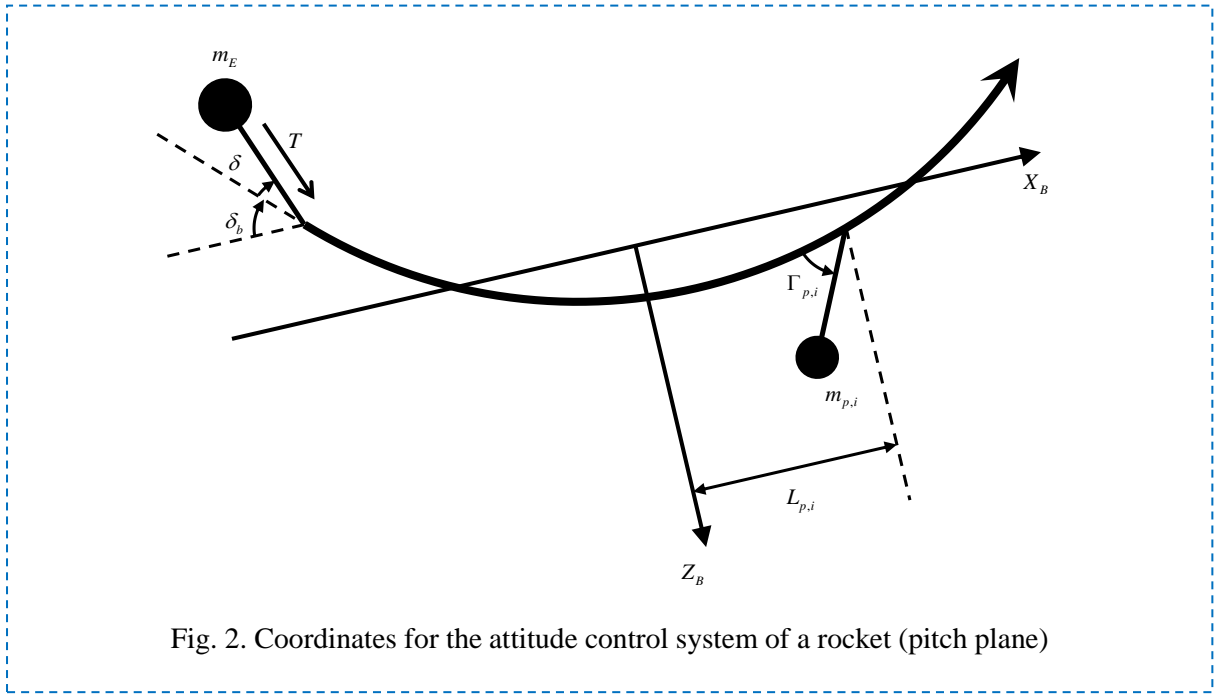

Fig. 2. Coordinates for the attitude control system of a rocket (pitch plane)

### 2.1.1. Rigid body motion with aerodynamics

A short period approximation model derived from the linearization of a nonlinear model about a nominal trajectory is used for the dynamic analysis of the rigid body motion with aerodynamics (Greensite, 1970). It is assumed that the rocket of interest is axisymmetric and the linearized equations of pitch motion are applicable to yaw motion. The linearized equations for the pitch motion of a rocket are expressed as follows:

$$\dot{\alpha} = \frac{\Delta F_z}{mU} \equiv Z = Z_\alpha \alpha + q + Z_\delta \delta + Z_S + Z_E \tag{1}$$

$$\dot{q} = \frac{\Delta M_y}{I_{yy}} \equiv M = M_\alpha \alpha + M_q q + M_\delta \delta + M_S + M_E \tag{2}$$

$$\dot{\theta} = q \tag{3}$$

where $\alpha$ is the angle of attack, $q$ is the pitch rate, $\theta$ is the pitch angle, $\delta$ is the TVC angle (in pitch direction) of a rocket. Note that these variables represent the small perturbations. The coefficients of $Z_\alpha$, $Z_\delta$, $M_\alpha$, $M_q$, and $M_\delta$ are given as follows:

$$Z_\alpha = \frac{\partial Z}{\partial \alpha} = -\frac{C_{N_\alpha} q_\infty S_{ref}}{mU} \tag{4}$$

$$Z_\delta = \frac{\partial Z}{\partial \delta} = \frac{T}{mU} \tag{5}$$

$$M_\alpha = \frac{\partial M}{\partial \alpha} = \frac{q_\infty C_{m_\alpha} S_{ref} l_{ref}}{I_{yy}} \tag{6}$$

$$M_q = \frac{\partial M}{\partial q} = \frac{C_{m_q} l_{ref}^2 S_{ref} \rho U}{4 I_{yy}} \tag{7}$$

$$M_\delta = \frac{\partial M}{\partial \delta} = \frac{(x_t - x_{cg}) \cdot T}{I_{yy}} \tag{8}$$

where $I_{yy}$ is the rotational inertia in $y$ axis, $U$ is the speed of the launch vehicle, $\rho$ is the air density, $q_\infty$ is the dynamic pressure ($= \rho U^2 / 2$), $S_{ref}$ and $l_{ref}$ are respectively the reference area and the reference length, $x_t$ and $x_{cg}$ are respectively the $x$ components of the thrust point and the center of gravity, and $T$ is the thrust. The derivatives of aerodynamic coefficients ($C_{N_\alpha}$, $C_{m_\alpha}$, and $C_{m_q}$) are defined as follows:

$$C_{N_\alpha} = \frac{\partial C_N}{\partial \alpha} \tag{9}$$

$$C_{m_\alpha} = \frac{\partial C_m}{\partial \alpha} \tag{10}$$

$$C_{m_q} = \frac{\partial C_m}{\partial (q l_{ref} / 2U)} \tag{11}$$

In addition, $Z_E / M_E$ are the TWD force/moment owing to the engine inertia and $Z_S / M_S$ are the force/moment resulting from the propellant sloshing, which are explained in detail in next subsections.

### 2.1.2. Structural bending of the rocket

Because of the slender shape of a rocket, body bending occurs in the rocket attitude motion, which can cause pitch/yaw instability. In particular, the flexible body rate due to the bending motion at the location of the gyro ($q_b$) is added to the pure pitch rate of the rigid body motion ($q$), which affects the stability of the attitude control system. It is known that a pair of the complex poles corresponding to the first bending mode are created close to the imaginary axis in the $s$ domain because of this additional measurement, and an attitude controller designed without considering these bending mode poles can lead to the instability of the system (Blakelock, 1991).

The bending mode can be described by applying Lagrange's equation to the rocket modeled as a system of lumped parameters (Greensite, 1970). The pitch rate input to the gyro sensor ($q_m$) is expressed as follows:

$$q_m = q + q_b = q + \sum_{i=1}^{N} \sigma_i(x_G)\eta_i \tag{12}$$

where $q_b$ is the local pitch rate caused by the bending motion, $\sigma_i(x_G)$ is the mode slope of the $i^{th}$ bending mode at the location of gyro ($x_G$), and $\eta_i$ is the generalized displacement of the $i^{th}$ bending mode expressed as the following second-order differential equation.

$$\ddot{\eta}_i + 2\zeta_{b,i}\omega_{b,i}\dot{\eta}_i + \omega_{b,i}^2\eta_i = -\frac{F_{G,i}}{M_{G,i}} \tag{13}$$

In this equation, $\omega_{b,i}$, $\zeta_{b,i}$, and $M_{G,i}$ are respectively the natural frequency, the damping ratio, and the generalized mass of the $i^{th}$ bending mode. In addition, the generalized force $F_{G,i}$ is calculated as follows:

$$F_{G,i} = \int_0^L F_z\phi_i(x)dx = \int_0^L (F_{Gz} + F_{Tz} + F_{Az} + F_{Sz} + F_{Ez})\phi_i(x)dx \tag{14}$$

where $\phi_i(x)$ is the mode shape of the $i^{th}$ bending mode. The force terms inside of the integral in Eq. (14) are the z-directional components of the gravitational force ($F_{Gz}$), thrust ($F_{Tz}$), aerodynamic force ($F_{Az}$), sloshing force ($F_{Sz}$), and engine's inertia force ($F_{Ez}$). First three of them are expressed as

$$F_{Gz} = mg\cos\theta \tag{15}$$

$$F_{Tz} = T\sin(\delta + \delta_b) \tag{16}$$

$$F_{Az} = -q_\infty S_{ref} C_{N_\alpha} (\alpha + \alpha_b) \tag{17}$$

where $\delta$ is the amount of TVC angle, and $\alpha_b$ and $\delta_b$ are respectively the angle of attack and TVC angle increments induced from the bending motion. The expressions for $F_{Ez}$ and $F_{Sz}$ are provided in next subsections.

### 2.1.3. TWD mode model

The effect of the inertial reaction force due to engine deflection, referred to the "Tail Wags Dog (TWD)" mode, is described as follows:

$$F_{Ez} = m_E \ell_E \ddot{\delta} \tag{18}$$

where $\ell_E$ is the distance from the pivot of the engine gimbal to the center of mass of the engine and $m_E$ is the mass of the engine. Then, the force and moment resulting from the engine's inertia in Eqs. (1) and (2) are expressed as follows:

$$Z_E = \frac{F_{Ez}}{mU} = \frac{m_E \ell_E}{mU} \ddot{\delta} \equiv Z_{\ddot{\delta}} \ddot{\delta} \tag{19}$$

$$M_E = \frac{I_E + m_E \ell_E \ell_C}{I_{yy}} \ddot{\delta} = M_{\ddot{\delta}} \ddot{\delta} \tag{20}$$

where $I_E$ is the moment of inertia of the engine, and $\ell_C$ is the distance from the center of mass of the rocket to the pivot of the engine gimbal.

### 2.1.4. Sloshing mode model

In a liquid propellant rocket, fuel sloshing produces forces and moments affecting the stability of its control system. The sloshing is modeled as a simple pendulum whose parameters can be derived from a partial differential equation describing the motion of liquid propellant in a cylindrical fuel tank with certain boundary conditions (Greensite, 1970). The sloshing motion of the propellant (fuel and oxidizer) in the $i^{th}$ tank is expressed as a differential equation on the pendulum angle $\Gamma_i$ as follows:

$$\ddot{\Gamma}_i + 2\zeta_{s,i}\omega_{s,i}\dot{\Gamma}_i + \omega_{s,i}^2\Gamma_i = -\frac{1}{L_{p,i}}\left[(\dot{W}-g_z) - Uq + Vp - \dot{q}(\ell_{p,i}-L_{p,i}) + \sum_{j=1}^{n_M}\ddot{\eta}_j\phi_j(\ell_{p,i})\right] \quad (21)$$

where $V$ and $W$ are the $y$ and $z$ components of the velocity, $p$ is the roll rate, $\ell_{p,i}$ is the distance between the center of mass of the rocket without propellant and the pivot of the pendulum, and $L_{p,i}$ is the length of the pendulum. The sloshing frequency ($\omega_{s,i}$) and the pendulum length ($L_{p,i}$) are expressed as follows:

$$\omega_{s,i} = \sqrt{\frac{a_x}{R_i}\xi_1\tanh\left(\frac{h_i}{R_i}\xi_1\right)} \quad (22)$$

$$L_{p,i} = \frac{a_x}{\omega_{s,i}^2} \quad (23)$$

where $a_x$ is the acceleration in $x$-direction, $\xi_1$ is a constant parameter (= 1.84), $R_i$ is the diameter of the tank, and $h_i$ is height of propellant, respectively. Then, the sloshing force $F_{Sz}$ in Eq. (14) is expressed as

$$\begin{aligned}
F_{Sz} &= \sum_{i=1}^{n_S} m_{p,i}L_{p,i}\left(\omega_{s,i}^2\Gamma_i + 2\zeta_{s,i}\omega_{s,i}\dot{\Gamma}_i\right) \\
&= \sum_{i=1}^{n_S} a_x m_{p,i}\left(\Gamma_i + \frac{2\zeta_{s,i}}{\omega_{s,i}}\dot{\Gamma}_i\right)
\end{aligned} \quad (24)$$

In this equation, $m_{p,i}$ is the $i^{th}$ slosh mass expressed as follows

$$m_{p,i} = \left(\frac{2R_i\tanh(h_i\xi_1/R_i)}{h_i\xi_1(\xi_1^2-1)}\right)m \quad (25)$$

where $m$ is total propellant mass. The expressions for $Z_S$ and $M_S$ in Eqs. (1) and (2) are presented as follows:

$$Z_S = \frac{F_{Sz}}{mU} = \sum_{i=1}^{n_S}\frac{a_x m_{p,i}}{mU}\left(\Gamma_i + \frac{2\zeta_{s,i}}{\omega_{s,i}}\dot{\Gamma}_i\right) \equiv \sum_{i=1}^{n_S}\left(Z_{\Gamma_i}\Gamma_i + Z_{\dot{\Gamma}_i}\dot{\Gamma}_i\right) \quad (26)$$

$$M_S = \sum_{i=1}^{n_S} \left[ -\frac{a_x m_{p_i} l_{p,i}}{I_{yy}} \left( \Gamma_i + \frac{2\zeta_{s,i}}{\omega_{s,i}} \dot{\Gamma}_i \right) \right] \equiv \sum_{i=1}^{n_S} \left( M_{\Gamma_i} \Gamma_i + M_{\dot{\Gamma}_i} \dot{\Gamma}_i \right) \tag{27}$$

## 2.2. TVC actuator and gyro model

The actuator dynamics of the TVC system can be modeled as a proper transfer function relating the command angle $\delta_c$ to the actual TVC angle $\delta$. In this study, we use the following 6th order transfer function.

$$\frac{\delta}{\delta_c} = \frac{\omega_{a1}^2}{s^2 + 2\zeta_{a1}\omega_{a1}s + \omega_{a1}^2} \cdot \frac{\omega_{a2}^2}{s^2 + 2\zeta_{a2}\omega_{a2}s + \omega_{a2}^2} \cdot \frac{\omega_{a3}^2}{\omega_{a4}^2} \frac{s^2 + 2\zeta_{a4}\omega_{a4}s + \omega_{a4}^2}{s^2 + 2\zeta_{a3}\omega_{a3}s + \omega_{a3}^2} \tag{28}$$

where the four natural frequency and damping ratio values are determined based on the characteristics of the actuator.

A gyro can also be modeled as a proper transfer function relating the measured pitch rate $q_m (= q + q_b)$ to the gyro output $q_g$. This study adopts the following 2nd order transfer function.

$$\frac{q_g(s)}{q_m(s)} = \frac{\omega_g^2}{s^2 + 2\zeta_g\omega_g s + \omega_g^2} \tag{29}$$

where $\omega_g$ and $\zeta_g$ are respectively the natural frequency and the damping ratio characterizing the gyro measurement.

## 2.3. Flight computer

A flight computer takes gyro measurements, conducts gyro error compensation, calculates the TVC angle command, and delivers the command to actuators through a DAC with ZOH. The gyro measurements are digitized by a V/F converter with sampling. For the analysis of the control system, the following third order transfer function is used to model the converter.

$$H_{VF}(s) = \frac{1 - b_1\tau_{VF}s + b_2\tau_{VF}^2 s^2 - b_3\tau_{VF}^3 s^3}{1 + a_1\tau_{VF}s + a_2\tau_{VF}^2 s^2 + a_3\tau_{VF}^3 s^3} \tag{30}$$

where $\tau_{VF}$ is the sampling period of the V/F and $a_1$, $a_2$, $a_3$, $b_1$, $b_2$, and $b_3$ are positive coefficients. The gyro error compensation and command angle calculation are modeled as a simple time delay system with time delay $\tau_D$. A second-order Padé approximation is used as follows:

$$H_D(s) = \frac{1 - \dfrac{\tau_D}{2}s + \dfrac{\tau_D^2}{12}s^2}{1 + \dfrac{\tau_D}{2}s + \dfrac{\tau_D^2}{12}s^2} \tag{31}$$

Finally, the DAC with ZOH are approximated as a 2$^{nd}$ order proper transfer function as follows.

$$H_{DZ}(s) = \frac{1}{1 + \dfrac{\tau_{DZ}}{2}s + \dfrac{\tau_{DZ}^2}{12}s^2} \tag{32}$$

where $\tau_{DZ}$ is a time constant associated with the ZOH.

## 3. Structure of attitude control law and bending filter

An attitude control law is required to make a rocket follow a pre-defined pitch/yaw rate, and a bending filter is necessary for stabilizing an overall attitude control system. Fig. 3 shows the structure of the attitude control law with a bending mode stabilization filter considered in this paper, and the details of the attitude control law and the bending filter are described in the following subsections.

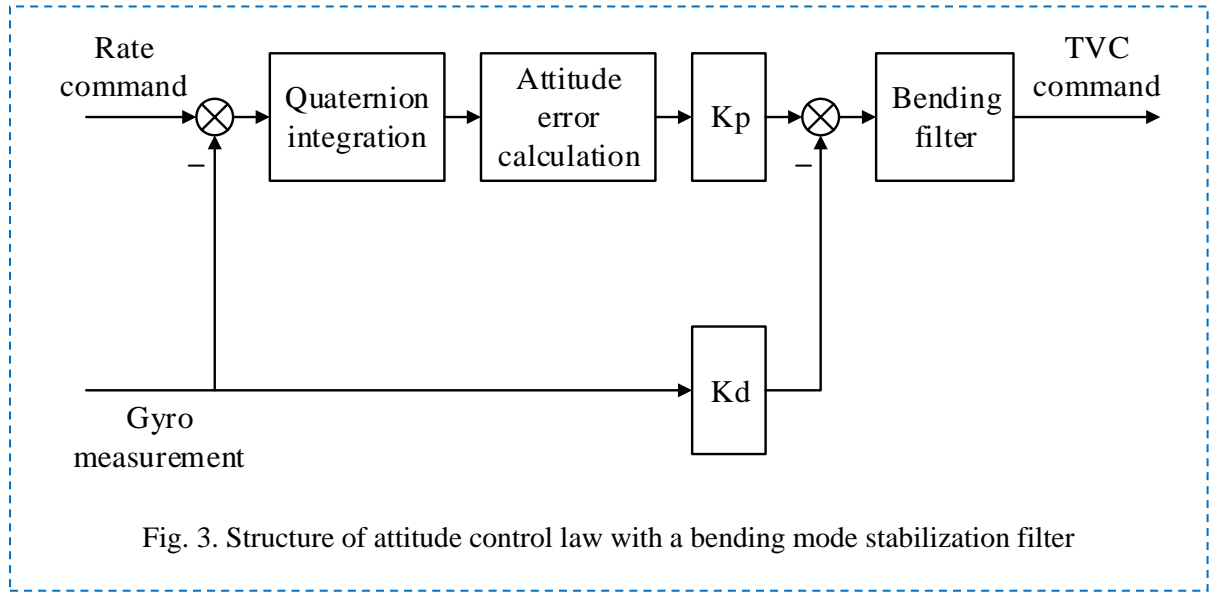

Fig. 3. Structure of attitude control law with a bending mode stabilization filter

### 3.1. Attitude control law

In this paper, we consider a simple proportional-derivative (PD) control law for the attitude control, and the TVC control input, which is filtered for bending mode stabilization, is expressed as follows.

$$\delta_c = -k_P \theta_e - k_D q_g \tag{33}$$

where $\theta_e \equiv \int_0^t q_e d\tau = \int_0^t (q_g - q_c) d\tau$.

When the PD control law is applied to the system dynamics (Eq. (2)) and only the rigid body motion is considered, the following characteristic equation can be obtained.

$$s^2 + M_\delta k_D s + M_\delta k_P = 0 \tag{34}$$

Thus, given desired values of the natural frequency and the damping ratio for the rigid body motion, one can obtain $k_P$ and $k_D$ from the following relationships derived from Eq. (34).

$$\omega^{RB} = \sqrt{M_\delta k_P} \tag{35}$$

$$\zeta^{RB} = \frac{1}{2} k_D \sqrt{\frac{M_\delta}{k_P}} \tag{36}$$

In general, $M_\delta$ varies during the flight time of a rocket, and as can be seen from Eqs. (35) and (36), it is impossible to keep the values of the natural frequency and damping ratio to be the desired ones with $k_P$ and $k_D$ fixed. Thus, a gain scheduling technique is adopted to design the values of $k_P$ and $k_D$ for each flight time segment.

## 3.2. Bending mode stabilization filter structure

The damping ratios of the structural bending modes are usually very small and large amplification around the natural frequency of each mode can occur. The bending mode, when not properly addressed, can cause the instability of the overall control system. A bending filter is designed and implemented to stabilize the system.

The bending filter shape the overall transfer function to secure a sufficient phase margin near the first bending mode (phase stabilization) and attenuate the magnitudes of higher-order bending modes (gain stabilization) while minimizing the low-frequency magnitude and phase. Five different filter structures ($F_2(s)$, $F_3(s)$, $F_4(s)$, $F_5(s)$, and $F_6(s)$) – with the orders from two to six – are considered in this study as follows:

$$F_2(s) = \frac{\dfrac{s^2}{\omega_z^2} + 2\dfrac{\zeta_z s}{\omega_z} + 1}{\dfrac{s^2}{\omega_p^2} + 2\dfrac{\zeta_p s}{\omega_p} + 1} \tag{37}$$

$$F_3(s) = \frac{\left(\dfrac{s^2}{\omega_z^2} + 2\dfrac{\zeta_z s}{\omega_z} + 1\right)\left(\dfrac{s}{z} + 1\right)}{\left(\dfrac{s^2}{\omega_p^2} + 2\dfrac{\zeta_p s}{\omega_p} + 1\right)\left(\dfrac{s}{p} + 1\right)} \tag{38}$$

$$F_4(s) = \frac{\left(\dfrac{s^2}{\omega_{z1}^2} + 2\dfrac{\zeta_{z1}s}{\omega_{z1}} + 1\right)\left(\dfrac{s^2}{\omega_{z2}^2} + 2\dfrac{\zeta_{z2}s}{\omega_{z2}} + 1\right)}{\left(\dfrac{s^2}{\omega_{p1}^2} + 2\dfrac{\zeta_{p1}s}{\omega_{p1}} + 1\right)\left(\dfrac{s^2}{\omega_{p2}^2} + 2\dfrac{\zeta_{p2}s}{\omega_{p2}} + 1\right)} \tag{39}$$

$$F_5(s) = \frac{\left(\dfrac{s^2}{\omega_{z1}^2} + 2\dfrac{\zeta_{z1}s}{\omega_{z1}} + 1\right)\left(\dfrac{s^2}{\omega_{z2}^2} + 2\dfrac{\zeta_{z2}s}{\omega_{z2}} + 1\right)\left(\dfrac{s}{z} + 1\right)}{\left(\dfrac{s^2}{\omega_{p1}^2} + 2\dfrac{\zeta_{p1}s}{\omega_{p1}} + 1\right)\left(\dfrac{s^2}{\omega_{p2}^2} + 2\dfrac{\zeta_{p2}s}{\omega_{p2}} + 1\right)\left(\dfrac{s}{p} + 1\right)} \tag{40}$$

$$F_6(s) = \frac{\left(\dfrac{s^2}{\omega_{z1}^2} + 2\dfrac{\zeta_{z1}s}{\omega_{z1}} + 1\right)\left(\dfrac{s^2}{\omega_{z2}^2} + 2\dfrac{\zeta_{z2}s}{\omega_{z2}} + 1\right)\left(\dfrac{s}{z_1} + 1\right)\left(\dfrac{s}{z_2} + 1\right)}{\left(\dfrac{s^2}{\omega_{p1}^2} + 2\dfrac{\zeta_{p1}s}{\omega_{p1}} + 1\right)\left(\dfrac{s^2}{\omega_{p2}^2} + 2\dfrac{\zeta_{p2}s}{\omega_{p2}} + 1\right)\left(\dfrac{s}{p_1} + 1\right)\left(\dfrac{s}{p_2} + 1\right)} \tag{41}$$

The schedules of gains ($k_P$ and $k_D$) of the control law introduced in Section 3.1 and the constant parameters of the bending filter in Section 3.2 are the design variables that should be determined by system engineers. Many published studies adopted a sequential approach; they first designs the gain schedules considering the rigid body motion (Eqs. (34)-(36)) with variations in the system parameters due to the large changes in inertia properties then design the bending filter (e.g., one of forms presented in Eqs. (37)-(41)) to satisfy the stability requirements such as the gain/phase margins. An optimization formulation for the integrated design of the gain schedules and bending filter is proposed in the next section.

# 4. Integrated design optimization formulation

The integrated gain schedule and bending filter design optimization problem for the attitude control of a rocket during its ascent flight can be formulated as the following nonlinear programming:

$$\min_{\mathbf{k}, \mathbf{x}} \left[ \max_{r \in \mathbf{R}} \left\| G(j\omega_1^B; k_P^r, k_D^r, \mathbf{x}) \right\| \right] \tag{42}$$

subject to

$$\left\| G(j\omega_2^B; k_P(t_i), k_D(t_i), \mathbf{x}) \right\| \leq P_{2,UL}^B \quad \forall i \in \mathbf{I} \tag{43}$$

$$\text{GM}(G(s; k_P(t_i), k_D(t_i), \mathbf{x})) \geq \text{GM}_R \quad \forall i \in \mathbf{I} \tag{44}$$

$$\text{PM}(G(s; k_P(t_i), k_D(t_i), \mathbf{x})) \geq \text{PM}_R \quad \forall i \in \mathbf{I} \tag{45}$$

$$\left\| F_k(j\infty; \mathbf{x}) \right\| \leq P^F \tag{46}$$

$$\omega_{RB}(k_P^i, k_D^i) \geq \omega_{RB,\min} \quad \forall i \in \mathbf{I} \tag{47}$$

$$\zeta_{RB}(k_P^i, k_D^i) \geq \zeta_{RB,\min} \quad \forall i \in \mathbf{I} \tag{48}$$

$$\mathbf{k} \geq \mathbf{0} \tag{49}$$

$$\mathbf{x}_{LB} \leq \mathbf{x} \leq \mathbf{x}_{UB} \tag{50}$$

The objective of the optimization, expressed in Eq. (42), is to minimize the worst case (out of the whole flight time) peak of the open-loop transfer function at the first bending mode frequency ($\omega_1^B$). In this equation, $G(s; k_P, k_D, \mathbf{x})$ is the open-loop transfer function of the pitch attitude control system from $q_c(s)$ to $q(s)$ with given proportional/differential gains ($k_P$ and $k_D$) and filter parameters ($\mathbf{x}$). The decision variable $\mathbf{k}$ denotes the vector containing the proportional and derivative gains at pre-determined times to define their schedules as follows:

$$\mathbf{k} = \left[ \cdots, k_P^r, k_D^r, \cdots \right], r \in \mathbf{R} = \{1, \cdots, R\} \tag{51}$$

The gain schedules $k_P(t)$ and $k_D(t)$ are obtained by linearly interpolating the elements of $\mathbf{k}$. The other decision variable vector $\mathbf{x}$ contains the structural bending filter parameters described in Eqs. (37) -(41). For example, if the 6$^{\text{th}}$ order bending filter $F_6(s)$ presented in Eq. (41) is used, $\mathbf{x}$ is expressed

as follows:

$$\mathbf{x} = [\omega_{z1}, \omega_{z2}, \zeta_{z1}, \zeta_{z2}, z_1, z_2, \omega_{p1}, \omega_{p2}, \zeta_{p1}, \zeta_{p2}, p_1, p_2] \tag{52}$$

Eq. (43) expresses the constraint that the second bending mode peak of the transfer function should not exceed its upper limit ($P_{2,UL}^B$). Eqs. (44) and (45) describe the constraints on the gain and phase margins of the control system at $t_i$, where $\mathrm{GM}(\cdot) / \mathrm{PM}(\cdot)$ are the gain margin / phase margin of a transfer function and $\mathrm{GM}_R / \mathrm{PM}_R$ are their minimum requirements, respectively. Eq. (46) accounts for the upper bound on the filter gain at infinite frequency. Eqs. (47) and (48) describes the constraints on the values of the natural frequency and the damping ratio for the rigid body motion (defined in Eqs. (35)-(36)) during the whole flight time. Finally, the bounds on the decision variable vectors are presented as inequality constraints in Eqs. (49) and (50).

# 5. Case Study

The integrated optimization framework for the gain schedules and bending filter was applied to the pitch/yaw control system design of Korea Sounding Rocket III (KSR-III), whose flight test was successfully completed in 2002. KSR-III is a single-stage liquid propellant sounding rocket launched in 2002. Table 1 summarizes the vehicle configuration and flight characteristics of KSR-III (Ahn et al., 2002; Ahn & Roh, 2012).

Table 1. Vehicle configuration and flight characteristics of KSR-III

| Vehicle configuration / Flight characteristic | Value/Property |
|---|---|
| Type | Sounding rocket (suborbital) |
| Nationality | Republic of Korea |
| Launch Date | Nov. 28, 2002 |
| Liftoff weight (kg) | 6,100 |
| Thrust (kN) | 129 |
| Burn time (s) | 59 |
| Flight time (s) | 230 |
| Max altitude (km) / range (km) | 43 / 79 |
| Max speed (m/s) / acceleration (g) | 850 / 4.1 |

Total 11 equally-spaced time nodes between 5 seconds and 55 seconds ($t_1$ = 5 sec., $t_2$ = 10 sec., …, $t_{11}$ = 55 sec.; $\mathbf{I}$ = {1, 2, …, 11}) were considered for checking the stability margins. Parameters for the models representing the dynamics of the rocket (presented in Section 2) are summarized in Table 2/Fig.4 (rigid body and aerodynamics), Table 3 (TVC actuator), and Table 4 (gyro and flight computer). Fig. 5 and Fig. 6 show the temporal variations of bending mode and sloshing frequencies. Note that for the sloshing mode, only the first mode is used because the forces and moments caused by the higher modes are very small (Blakelock, 1991).

Table 2. Parameter ranges for rigid body / aerodynamic models (stability derivatives)

| Parameter | $Z_\alpha$ ($s^{-1}$) | $M_\alpha$ ($s^{-2}$) | $Z_\delta$ ($s^{-1}$) | $M_\delta$ ($s^{-2}$) |
|---|---|---|---|---|
| Maximum value | -0.1 | -1.0 | 0.38 | 15.5 |
| Minimum value | -0.55 | -40 | 0.05 | 8.1 |

Table 3. Parameter values for TVC actuator model

| Parameter | $\omega_{a1}$ (Hz) | $\omega_{a2}$ (Hz) | $\omega_{a3}$ (Hz) | $\omega_{a4}$ (Hz) | $\zeta_{a1}$ (-) | $\zeta_{a2}$ (-) | $\zeta_{a3}$ (-) | $\zeta_{a4}$ (-) |
|---|---|---|---|---|---|---|---|---|
| Value | 5.8 | 7.2 | 52.6 | 15.6 | 0.87 | 0.49 | 0.15 | 0.042 |

Table 4. Parameter values for DTG / flight computer

| Model | DTG | | Flight computer | | |
|---|---|---|---|---|---|
| Parameter | $\omega_g$ (Hz) | $\zeta_g$ (-) | $\tau_{VF}$ (s) | $\tau_D$ (s) | $\tau_{DZ}$ (s) |
| Value | 60 | 0.65 | 0.01 | 0.003 | 0.01 |

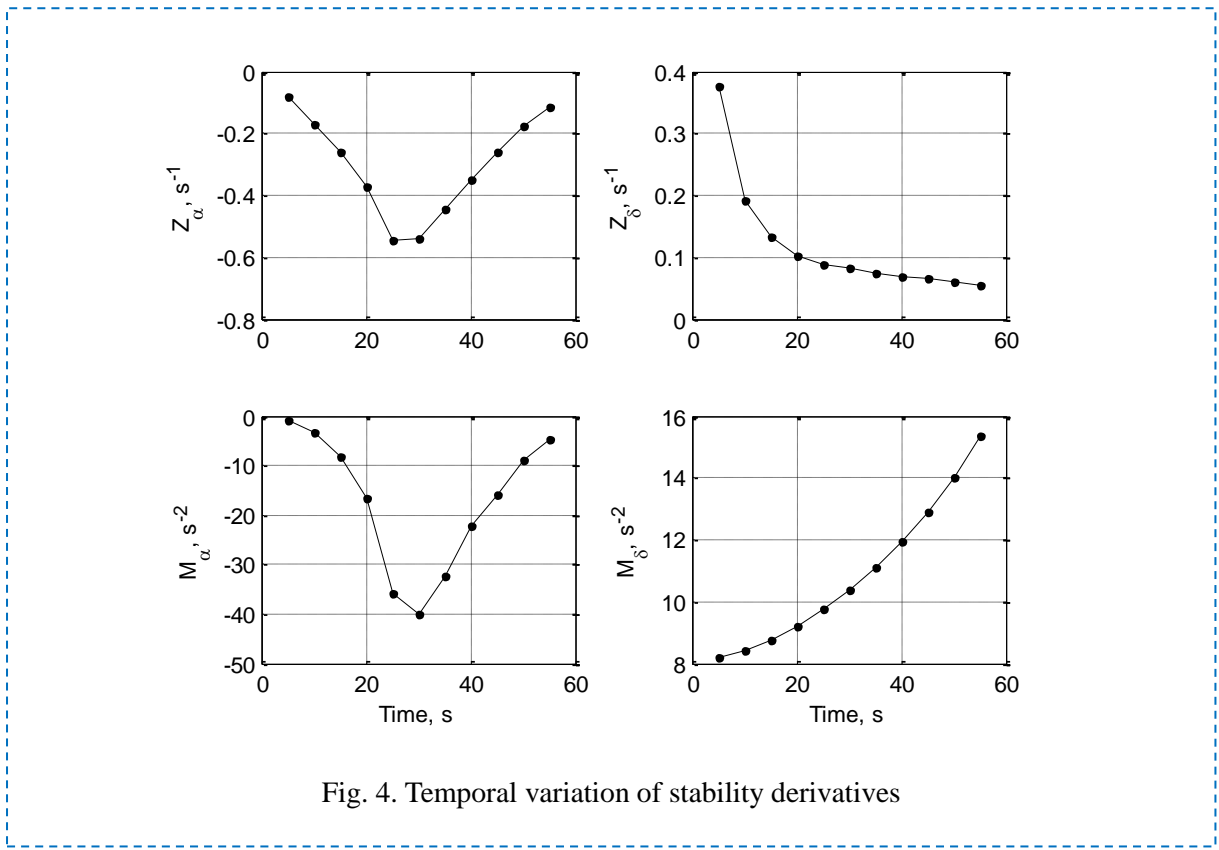

Fig. 4. Temporal variation of stability derivatives

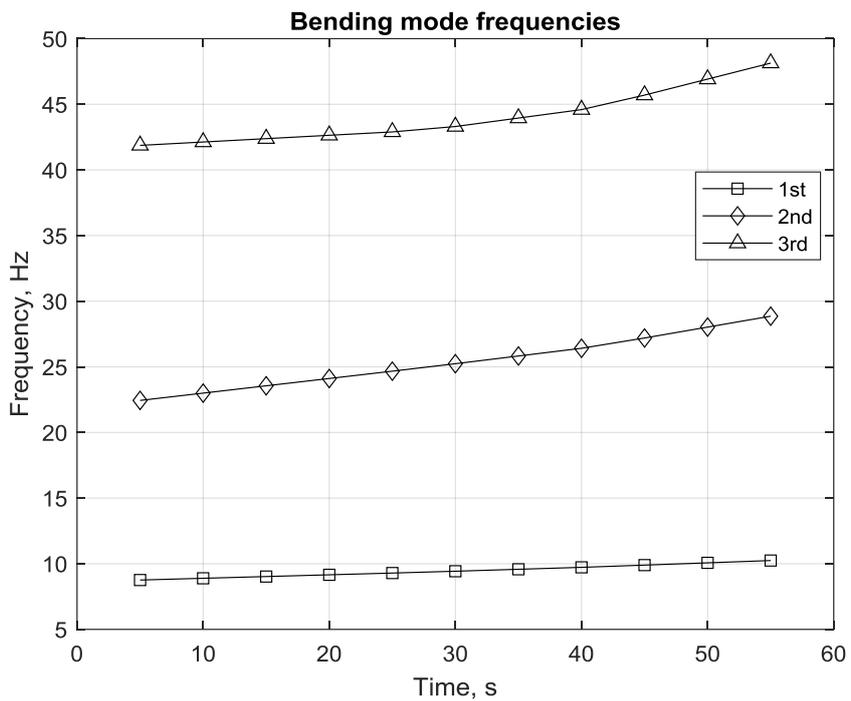

Fig. 5. First three bending mode frequencies

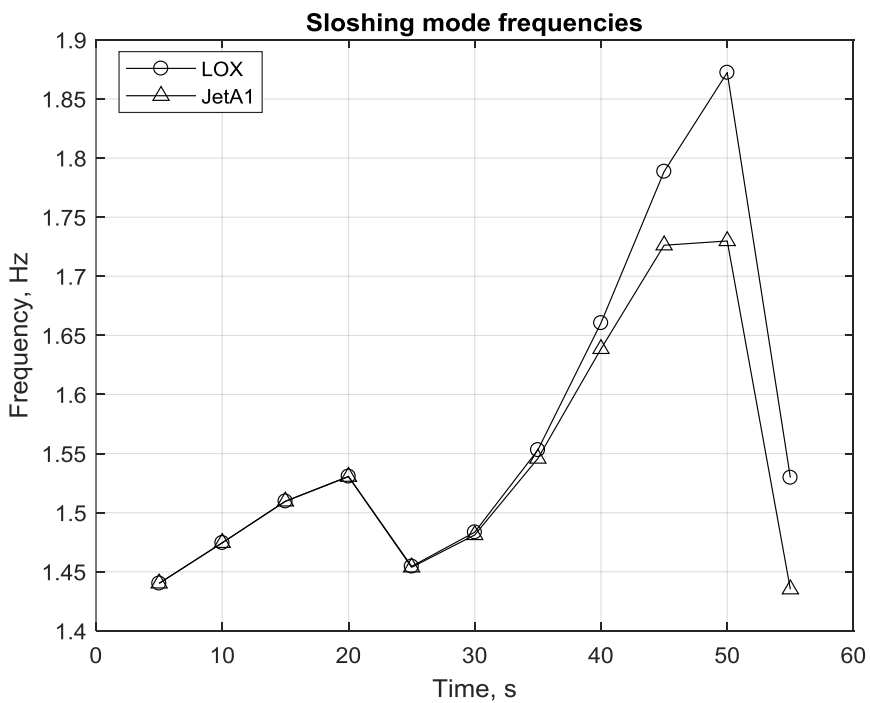

Fig. 6. Sloshing mode frequencies for LOX and JetA1

## 5.1. Reference case: two-phase design of gain schedules and bending filter

The performance of the gain schedules and bending filter obtained using the proposed procedure is compared with the sequential (two-phase) design results presented by Ahn et al. (2002), which were implemented for the actual flight test of KSR-III. In the first phase of the study, the gain schedules that vary linearly with respect to the flight time were designed. That is, the intercepts and slopes of the gain schedules were determined in such a way that the natural frequency and the damping ratio for the rigid body motion, $\omega_i^{RB}$ and $\zeta_i^{RB}$, are guaranteed to be no smaller than 0.5 Hz and 0.7, respectively. In the second phase, the bending mode stabilization filters described in Eqs. (37)-(41), were designed so that the worst-case (largest) first bending mode peak magnitude during the flight is minimized. The gain profiles and the Bode plot of the bending filters obtained by the two-phase design are shown in Fig. 7 and Fig. 8. The minimized worst-case first bending mode peak magnitudes associated with the designed filters of various orders are presented in Fig. 9 (diamond markers with a dotted line).

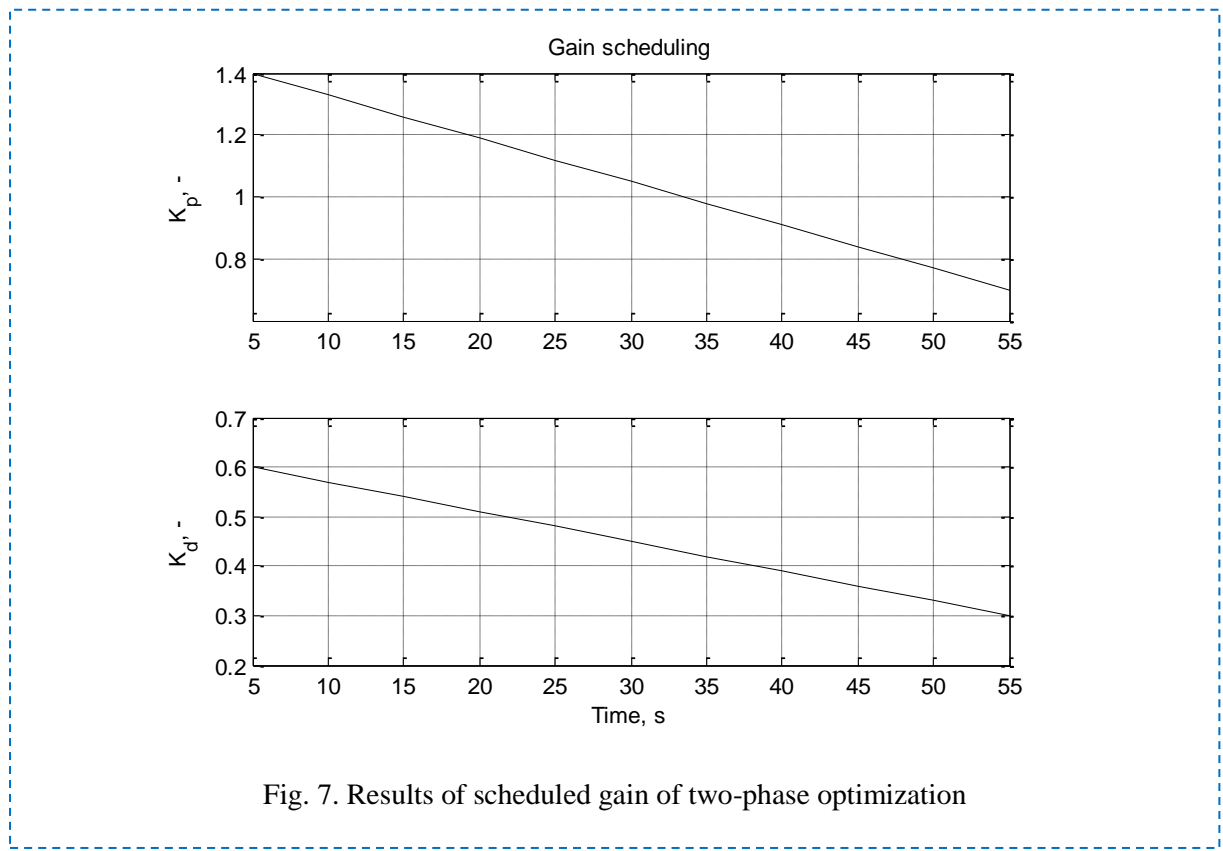

Fig. 7. Results of scheduled gain of two-phase optimization

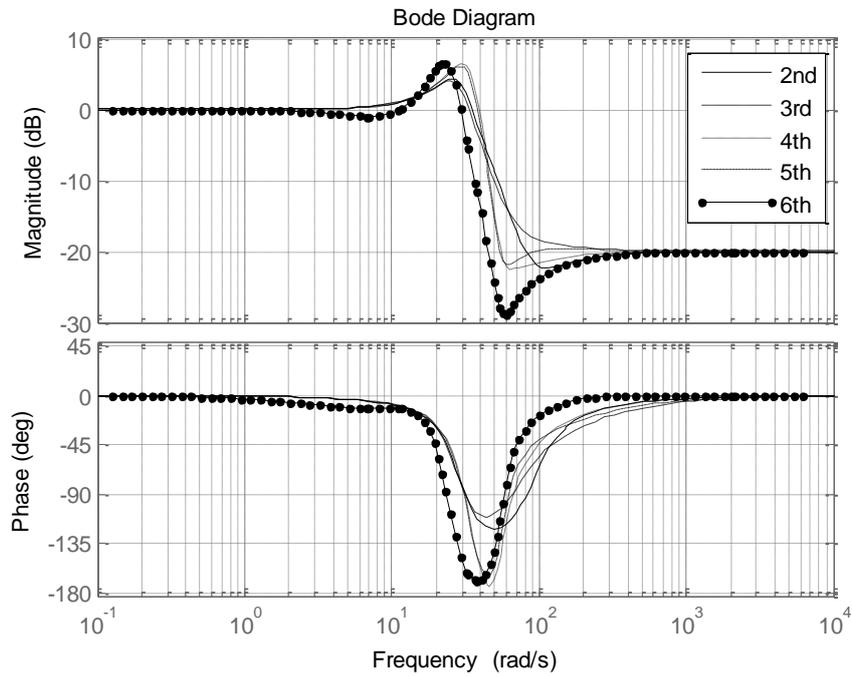

Fig. 8. Bode plot of optimized bending filters of two-phase optimization

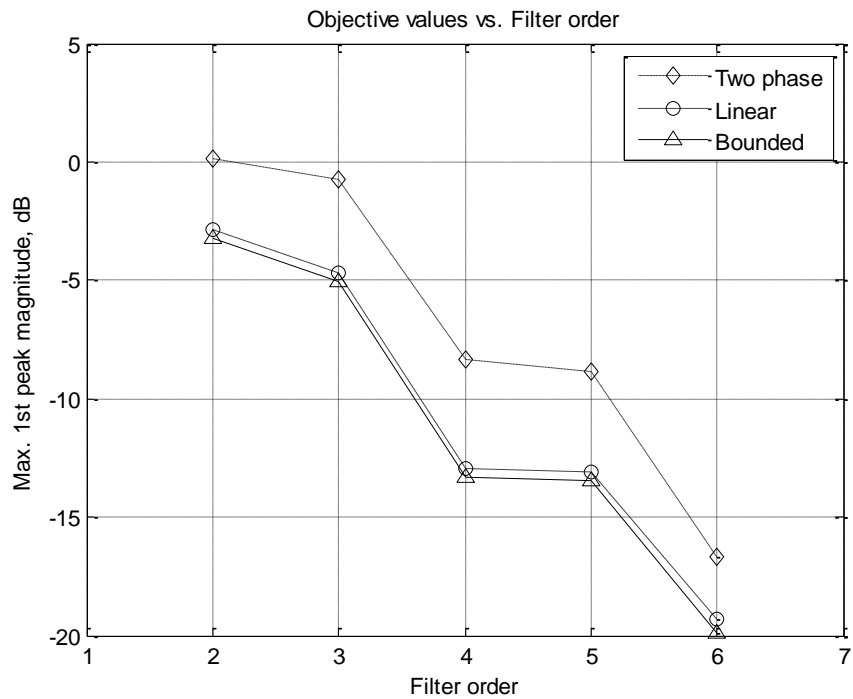

Fig. 9. Minimized maxima of the first bending mode peak magnitudes for the filers with various orders

## 5.2. Integrated design optimization of gain schedules and bending filter

The design of gain schedules and bending filter using the proposed integrated approach was conducted and its effectiveness was demonstrated through the comparison with the results of the reference case presented in Section 5.1. Two subcases with different gain schedule shapes (and associated decision variable $\mathbf{k}$) were considered. Subcase 1 assumes that the schedules of the proportional and differential gains ($k_P(t)$ and $k_D(t)$) are determined by the gain values at $t_1$ and $t_T$, which involves only two design variables for a gain schedule ($\mathbf{k} = [k_P^1, k_D^1, k_P^2, k_D^2]$). In Subcase 2, the values of proportional/differential gains at every $t_i$ ($i = 1, \ldots, T$) are selected as design variables ($\mathbf{k} = [k_P^1, k_D^1, \cdots, k_P^T, k_D^T]$), which allows for higher degree of freedom for improved performance of the attitude control system.

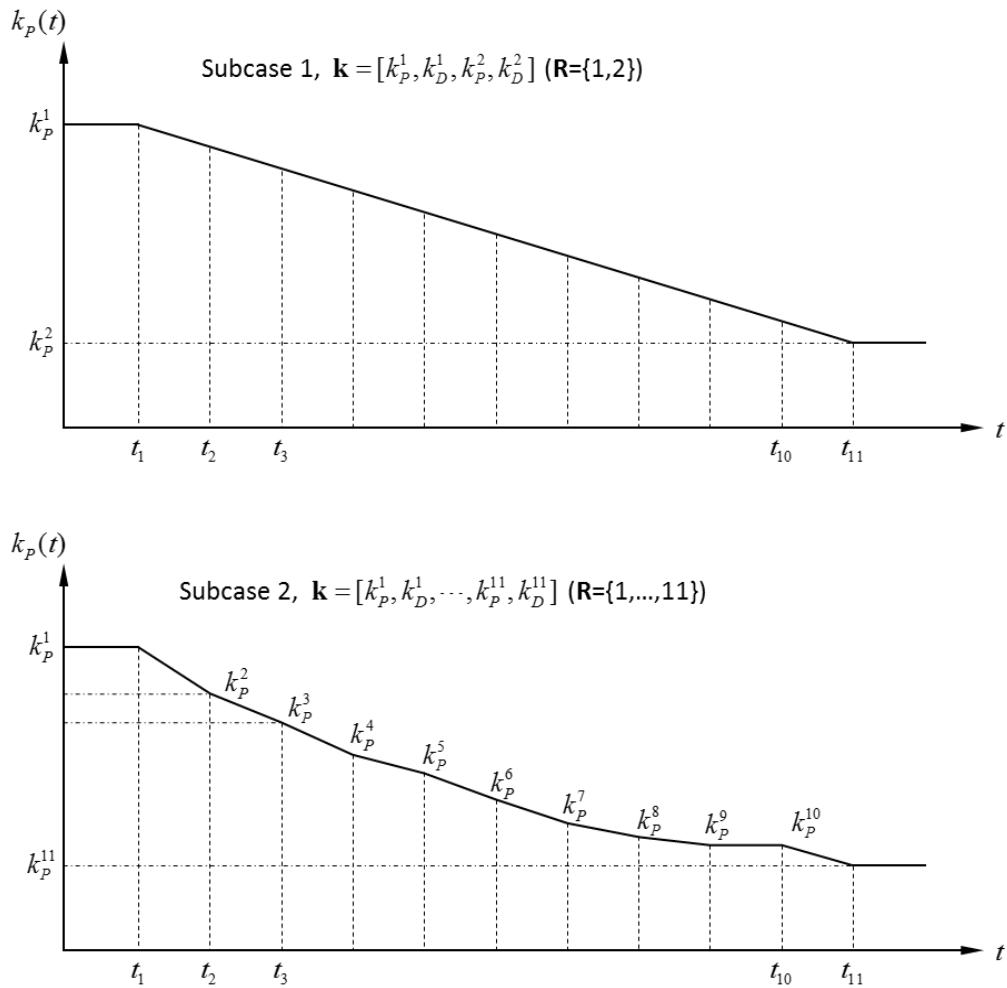

Fig. 10. Gain schedules and associated design variables for two subcases

The objective functions (the worst-case first bending mode peak during the flight) of the integrated optimization problem for various filter structures are presented in Fig. 9 (circle markers with a dashed line). Compared with the reference design, Subcase 1 can achieve the reduction of 3 to 5 dB in worst-case peak of the first bending mode for all filter structures. The gain schedule profiles and the Bode plots of the optimal bending filters for Subcase 1 are shown in Figs. 11 and 12, respectively. The comparison between Fig. 7 and Fig. 10 reveals that the optimal gain values obtained by the integrated design approach are generally smaller than reference gain values, which are determined by the sequential (two-phase) design procedure for the gain schedules and bending filter.

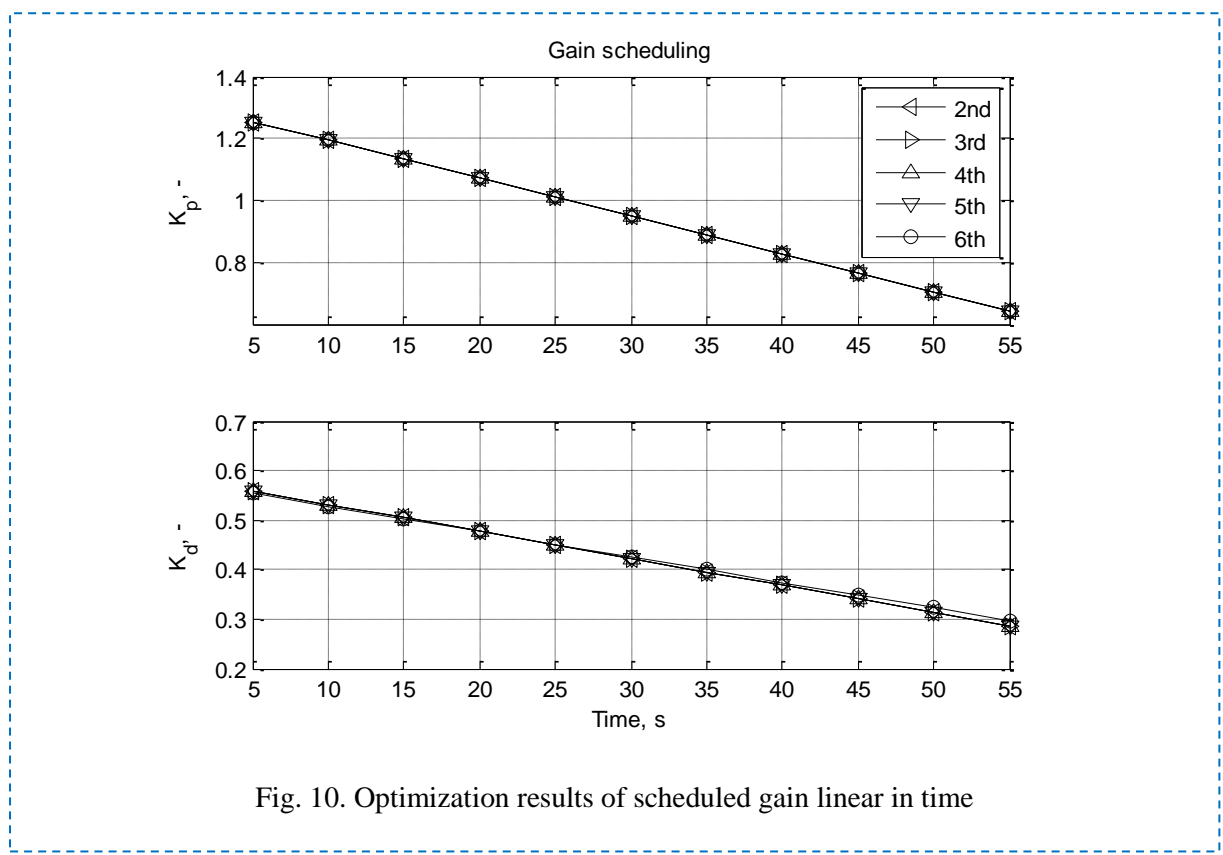

Fig. 10. Optimization results of scheduled gain linear in time

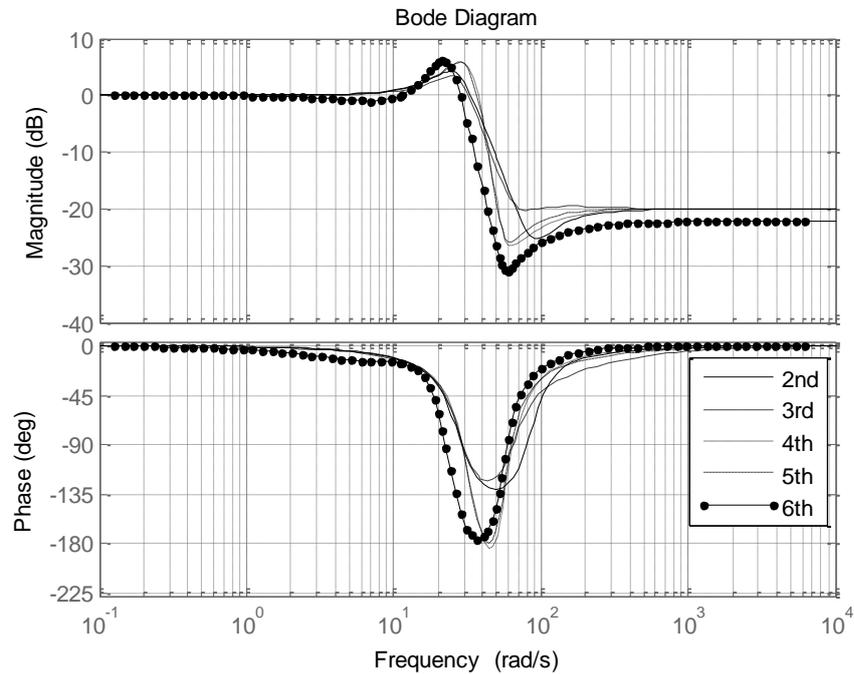

Fig. 12. Bode plot of optimized bending filters with ACS PD gains linear in time

The design results for Subcase 2 – the objective function, the gain schedules, and the Bode plot of the bending filters – are shown in Fig. 9 (triangular markers with a solid line), Fig. 13, and Fig. 14, respectively. As can be seen in Fig. 9, the objective function associated with the optimal design for Subcase 2 is the lowest for all filter structures. The amounts of reduction compared to Subcase 1 range between 0.33 to 0.57 dB, which is relatively insignificant considering the increase in the number of decision variables for determining the gain schedules (4 for Subcase 1 versus 22 for Subcase 2). The comparison indicates that, for this case, the small amount of improvement in the performance of the integrated gain schedule and filter design does not well justify the increased complexity for the optimal design problem for Subcase 2.

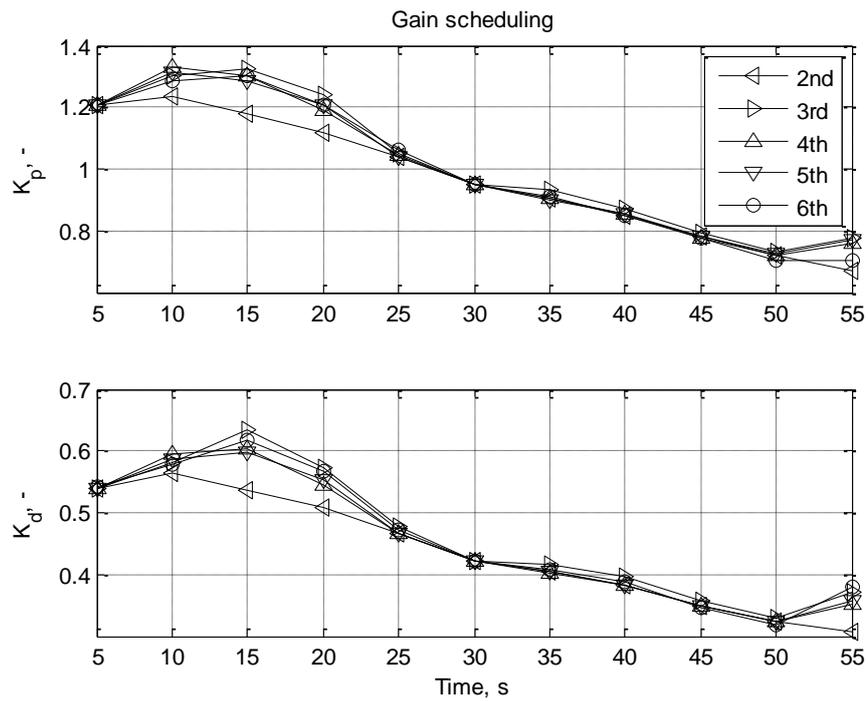

Fig. 13. Optimization results of scheduled gain bounded

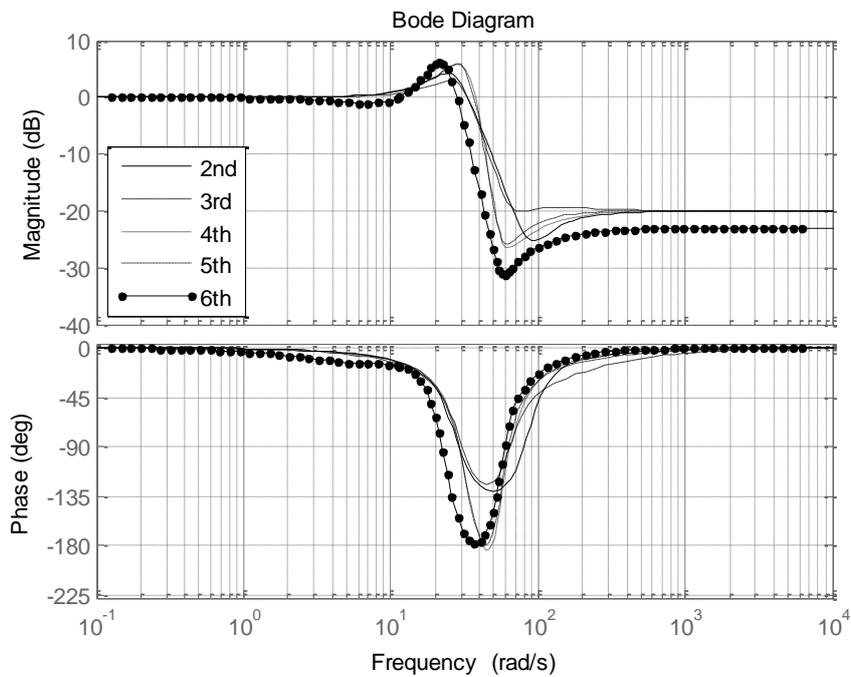

Fig. 14. Bode plot of optimized bending filters with ACS PD gains bounded

## 6. Conclusions

An optimization framework for the integrated design of gain schedules and bending filter for the pitch/yaw control of a rocket is proposed in this study. The proposed framework was applied to the case study on the longitudinal control system design of Korea Sounding Rocket III (KSR-III), which had flown in 2002. The performance of the obtained control system was compared with that of the sequential design of the gain schedules and bending filter, which had been implemented for the actual flight test of KSR-III. The comparison result indicates that the proposed procedure outperforms the sequential approach in terms of attenuation of the first bending mode peak of the open-loop loop transfer function by 3 to 5 dB, which demonstrated its applicability for control system design of rockets and launch vehicles.

## References

Ahn, J. & Roh, W. (2012). Noniterative instantaneous impact point prediction algorithm for launch operations. Journal of Guidance, Control, and Dynamics 35(2), 645-648. doi:10.2514/1.56395.

Ahn, J., Roh, W.-R., Cho, H.-C., & Park, J.-J. (2002). Control system modeling and optimal bending filter design for KSR-III first stage. *Journal of the Korean Society for Aeronautical and Space Sciences*, 30(7), 113-122. doi:10.5139/JKSAS.2002.30.7.113.

Blakelock, J. H. (1991). *Automatic Control of Aircraft and Missiles*. New York, NY: John Wiley & Sons, Inc.

Choi, H.-D. & Bang, H. (2000). An adaptive control approach to the attitude control of a flexible rocket. *Control Engineering Practice*, 8(9), 1003-1010. doi: 10.1016/S0967-0661(00)00032-0.

Clement, B., Duc, G., & Mauffrey, S. (2005). Aerospace Launch Vehicle Control: a Gain scheduling Approach. *Control Engineering Practice*, 13(3), 333-347. doi:10.1016/j.conengprac.2003.12.013.

Duraffourg, E., Burlion, L., & Ahmed-Ali, T. (2013). Longitudinal modeling and preliminary control of a nonlinear flexible launch vehicle. *11th IFAC International Workshop on Adaptation and Learning in Control and Signal Processing*, Caen.

Frosch, J. A. & Vallely, D. P. (1967). Saturn AS-501/S-IC flight control system design. *Journal of Spacecraft and Rocket* 4(8), 1003-1009.

Greensite, A. L. (1970). *Analysis and Design of Space Vehicle Flight Control Systems*. New York, NY: Spartan Books.

Haeussermann, W. (1970). Description and performance of the Saturn launch vehicle's navigation, guidance, and control system. NASA-TN-D-5869.

Jang, J.-W., Hall, R., & Bedrossian, N. (2008). Ares-I bending filter design using a constrained optimization approach. *AIAA Guidance, Navigation and Control Conference and Exhibit*, Honolulu. doi: 10.2514/6.2008-6289.

Kamath, A., Menon, P., Ganet-Schoeller, M., Maurice, G., Bennani, S., & Bates, D. (2012) Robust safety margin assessment and constrained worst-case analysis of a launcher vehicle. *7th IFAC Symposium on Robust Control Design*. Aalborg. doi: 10.3182/20120620-3-DK-2025.00136.

Mori, H. (1999). Control system design of flexible-body launch vehicles. *Control Engineering Practice*, 7, 1163-1175. doi:10.1016/S0967-0661(99)00087-8.

Morita, Y. & Kawaguchi, J. (2001). Attitude control design of the M-V rocket. Philosophical Transactions of the Royal Society A, 359(1788), 2287-2303. doi:10.1098/rsta.2001.0887.

Oh, C.-S., Bang, H., & Park, C.-S. (2008). Attitude control of a flexible launch vehicle using an adaptive notch filter: Ground experiment. Control Engineering Practice 16(1), 30-42. doi:10.1016/j.conengprac.2007.03.006.

Rao, K., Dhekane, M., Lalithambika, V., & Brinda, V. (2014). Application of optimization technique for compensator design of launcher attitude control. *Third International Conference on Advances in Control and Optimization of Dynamical Systems*, Kanpur. doi: 10.3182/20140313-3-IN-3024.00080.

Saunois, P. (2009). Comparative analysis of architectures for the control loop of launch vehicles during atmospheric flight. *Aerospace Science and Technology*, 13(2-3), 150-156. doi:10.1016/j.ast.2008.08.003.

Schleich, W. T. (1982). The Space shuttle ascent guidance and control. *AIAA Guidance and Control*

*Conference*. San Diego. doi:10.2514/6.1982-1497.

Simmons, J. C. & Hall, B. M. (1973). Digital control system development for the Delta launch

vehicle. *AIAA Guidance and Control Conference*. doi:10.2514/6.1973-847.